\documentclass[11pt]{article}
\usepackage{graphicx}
\usepackage{color}
\usepackage{amsmath}
\usepackage{amssymb}
\usepackage{amscd}
\usepackage{bbm}
\newcommand{\R}{\mathbb{R}}
\newcommand{\inr}[1]{\bigl< #1 \bigr>}

\newcommand{\E}{\mathbb{E}}

\newcommand{\eps}{\varepsilon}

\newtheorem{Theorem}{Theorem}[section]
\newtheorem{Lemma}[Theorem]{Lemma}

\newtheorem{Corollary}[Theorem]{Corollary}

\newtheorem{Question}[Theorem]{Question}
\numberwithin{equation}{section}

\def \proof {\noindent {\bf Proof.}\ \ }

\def \endproof
{{\mbox{}\nolinebreak\hfill\rule{2mm}{2mm}\par\medbreak}}

\begin{document}
\title{A remark on the diameter of random sections of convex bodies}
\author{Shahar Mendelson${}^{1}$}

\footnotetext[1]{Department of Mathematics, Technion, I.I.T, Haifa, 32000, Israel. Email: shahar@tx.technion.ac.il. Partially supported by the Mathematical Sciences Institute -- The Australian National University and by the Israel Science Foundation grant 900/10.}
\maketitle

\begin{abstract}
We obtain a new upper estimate on the Euclidean diameter of the intersection of the kernel of a random matrix with iid rows with a given convex body. The proof is based on a small-ball argument rather than on concentration and thus the estimate holds for relatively general matrix ensembles.
\end{abstract}

\section{Introduction}
In this note we revisit the following problem.

Let $\mu$ be an isotropic measure on $\R^n$, and by `isotropic' we mean a symmetric measure that satisfies
$$
\int_{\R^n} \inr{x,t}^2 d\mu(t) = \|x\|_{\ell_2^n}^2 \ {\rm \ for \ every}  \ x \in \R^n.
$$
Given a random vector $X$ distributed according to $\mu$ and for $X_1,...,X_k$ that are independent copies of $X$, let $\Gamma$ be the random matrix $k^{-1/2} \sum_{i=1}^k \inr{X_i,\cdot}e_i$.

\begin{Question}
If $T \subset \R^n$ is a convex body (that is, a convex, centrally-symmetric set with a non-empty interior), what is the typical diameter of $T \cap {\rm ker}(\Gamma)$?
\end{Question}

The origin of this problem was the study of the geometry of convex bodies, and in particular, Milman's low-$M^*$ estimate \cite{Mi} and subsequent estimates on the Gelfand widths of convex bodies, due to Pajor and Tomczak-Jaegermann \cite{PT1,PT2}.

The focus of the original question had been the existence of a section of $T$ of codimension $k$ and of a small Euclidean diameter, and was established by estimating $\E {\rm diam}(T \cap E)$ from above, relative to the uniform measure on the Grassmann manifold $G_{n-k,n}$.

In recent years, more emphasis has been put on other choices of measures on the Grassmann manifold, for example, using the distribution generated by kernels of matrices selected from some random ensemble -- like $\Gamma=k^{-1/2}\sum_{i=1}^k \inr{X_i,\cdot}e_i$ defined above.

The standard way of estimating $\E {\rm diam}(T \cap {\rm ker}(\Gamma))$ for such matrix ensembles is based on the quadratic empirical processes indexed by linear forms associated with $T$.

It is straightforward to show (see, for example, the discussion in \cite{MPT}) that  given $r>0$, if
\begin{equation} \label{eq:quadratic}
\sup_{x \in T \cap r S^{n-1}} \left|\frac{1}{k}\sum_{i=1}^k \inr{X_i,x}^2 - \E\inr{X,x}^2 \right| \leq \frac{r^2}{2},
\end{equation}
one has
$$
\frac{1}{2} \|x\|_{\ell_2^n}^2 \leq \|\Gamma x\|_{\ell_2^k}^2 \leq \frac{3}{2}\|x\|_{\ell_2^n}^2
$$
for every $x \in T$ of $\ell_2^n$ norm larger than $r$. Hence, on the event given by \eqref{eq:quadratic}, ${\rm diam}(T \cap {\rm ker}(\Gamma)) \leq r$.

Setting $r_0(k,\delta)$ to be the smallest for which
$$
Pr \left(\sup_{x \in T \cap rS^{n-1}} \left|\frac{1}{k}\sum_{i=1}^k \inr{X_i,x}^2 - \E\inr{X,x}^2 \right| \leq \frac{r^2}{2}\right) \geq 1-\delta,
$$
it follows that with probability at least $1-\delta$,
$$
{\rm diam}(T \cap {\rm ker}(\Gamma)) \leq r_0,
$$
and a similar argument may be used to control $\E {\rm diam}(T \cap {\rm ker}(\Gamma))$.

Unfortunately, estimating the quadratic empirical process is a difficult task. In fact, one has a satisfactory estimate that holds for every convex body $T \subset \R^n$ only for measures that are subgaussian or unconditional log-concave.

\begin{Theorem} \label{thm:MPT} \cite{Men-dvo}
There exist absolute constants $c_1,c_2$ and $c_3$ for which the following holds. Let $\mu$ be an isotropic, $L$-subgaussian measure (and, in particular, for every $x \in \R^n$, $\|\inr{X,x}\|_{\psi_2(\mu)} \leq L \|x\|_{\ell_2^n}$).

Let $T \subset \R^n$ and set $d_T = \sup_{t \in T} \|t\|_{\ell_2^n}$.
For $u \geq c_1$, with probability at least
$$
1-2\exp(-c_2u^2(\E\|G\|_{T^\circ} /d_T)^2),
$$
\begin{equation*}
\sup_{x \in T} \left|\sum_{i=1}^k (\inr{X_i,x}^2 - \E|\inr{X,x}|^2) \right|
\leq c_3L^2 u^2 \left((\E\|G\|_{T^\circ})^2 + u\sqrt{k}d_T\E\|G\|_{T^\circ}\right),
\end{equation*}
where $G=(g_1,...,g_n)$ is the standard gaussian vector in $\R^n$ and $\E\|G\|_{T^\circ} = \E \sup_{t \in T} |\inr{G,t}|$.
\end{Theorem}
A version of Theorem \ref{thm:MPT} has been established in \cite{MPT} when $T \subset S^{n-1}$ and with a weaker probability estimate.

Theorem \ref{thm:MPT} follows from a general bound on the quadratic empirical process that is based on a global complexity parameter of the indexing set \cite{MenPao}, and that will not be defined here. Thanks to Talagrand's Majorizing Measures Theorem (see the book \cite{Tal-book} for a detailed survey on this topic), this complexity parameter is upper bounded by $\sim \E\|G\|_{T^\circ}$ in the subgaussian case, thus leading to Theorem \ref{thm:MPT}. However, in other cases, controlling it is nontrivial.

One other case in which the global complexity may be upper bounded using a mean-width of $T$, is when $X$ is isotropic, unconditional and log-concave. Using the Bobkov-Nazarov Theorem \cite{Bob-Naz}, $X$ is dominated by $Y=(y_1,...,y_n)$, a vector with independent, standard, exponential coordinates. One may show \cite{MenPao} that with high probability,
\begin{equation} \label{eq:exp}
\sup_{x \in T} \left|\sum_{i=1}^k (\inr{X_i,x}^2 - \E|\inr{X,x}|^2) \right|
\lesssim  (\E\|Y\|_{T^\circ})^2 + \sqrt{k}d_T\E\|Y\|_{T^\circ}.
\end{equation}
The proof of \eqref{eq:exp} is based on two additional observations. First, that when $X$ is isotropic, unconditional and log-concave, the global complexity parameter of $T$ may be bounded using a mixture of Talagrand's $\gamma_\alpha$ functionals, and second, that this mixture is equivalent to $\E\|Y\|_{T^\circ}$ \cite{Tal-AJM}.

Additional bounds on the quadratic process are known for more general measures, but only for very specific choices of sets $T$. The most important example is when $T$ is the Euclidean ball, and the quadratic empirical process may be used to obtain a Bai-Yin type estimate on the largest and smallest singular values of $\Gamma$ \cite{MenPao,MenPao1}.

At this point, it should be noted that \eqref{eq:quadratic} is a much stronger statement than what is actually needed to bound the diameter of $T \cap {\rm ker}(\Gamma)$. Clearly, any sort of a positive {\it lower} bound on
\begin{equation} \label{eq:inf}
\inf_{x \in T \cap rS^{n-1}} \|\Gamma x\|_{\ell_2^k}
\end{equation}
would suffice -- rather than the `almost isometric', two-sided bound that follows from bounds on the quadratic process.

Here, we will show that \eqref{eq:inf} holds for rather general matrix ensembles.

\begin{Theorem} \label{thm:main}
Let $X$ be an isotropic vector on $\R^n$ and assume that linear forms satisfy the following small-ball condition: that there is some $\lambda>0$ for which
$$
Pr(|\inr{x,X}| \geq \lambda\|x\|_{\ell_2^n} ) \geq 99/100 \ \ \ {\rm for \ every \ } x \in \R^n.
$$
Then, there exist a constant $c$ that depends only on $\lambda$, for which, with probability at least $3/4$,
$$
{\rm diam}(T \cap {\rm ker}(\Gamma)) \leq \frac{c}{\sqrt{k}} \cdot  \max\left\{\E\|G\|_{T^\circ}, \E\|k^{-1/2}\sum_{i=1}^k X_i\|_{T^\circ} \right\}.
$$
\end{Theorem}

Theorem \ref{thm:main} can be improved and extended in various ways.

First of all, the `correct' upper estimate on the diameter should be based on a fixed point condition defined using the norms $\| \ \|_{(T \cap rB_2^n)^\circ}$ rather than the norm $\| \ \|_{T^\circ}$. Also, the constant probability estimate of $3/4$ may be improved significantly to $1-2\exp(-ck)$ with a slightly more involved proof (see \cite{KM} for a similar argument). We will formulate, without proof, a more general version of Theorem \ref{thm:main} at the end of the note.
\vskip0.5cm
\noindent{\bf Examples.}

\begin{description}
\item{1.} If $X$ is an isotropic $L$-subgaussian vector, it is standard to verify that $k^{-1/2}\sum_{i=1}^k X_i$ is isotropic and $cL$-subgaussian for a suitable absolute constant $c$. Therefore,
$$
\E\|k^{-1/2} \sum_{i=1}^k X_i\|_{T^\circ} \leq c_1L \E\|G\|_{T^\circ},
$$
and by Theorem \ref{thm:main}, with probability at least $3/4$,
$$
{\rm diam}(T \cap {\rm ker}(\Gamma)) \leq c_1(\lambda,L) \frac{\E \|G\|_{T^\circ}}{\sqrt{k}}.
$$
This coincides with the estimate from \cite{MPT} (up to the `localization' mentioned above) and with the classical result of \cite{PT1} when $X$ is the standard gaussian vector.

\item{2.}
 If $X$ is an isotropic, unconditional, log-concave measure then so is $Z=k^{-1/2} \sum_{i=1}^k X_i$. By the Bobkov-Nazarov Theorem \cite{Bob-Naz}, both $Z$ and $G$ are strongly dominated by $Y$, the random vector with independent, standard exponential coordinates. Therefore, by Theorem \ref{thm:main}, with probability at least $3/4$,
$$
{\rm diam}(T \cap {\rm ker}(\Gamma)) \leq c_2(\lambda)  \frac{\E\|Y\|_{T^\circ}}{\sqrt{k}}.
$$

\item{3.} Theorem \ref{thm:main} leads to a `heavy tails' result in some cases. Since $X$ is symmetric, $\sum_{i=1}^k X_i$ has the same distribution as $\sum_{i=1}^k \eps_i X_i$, where $(\eps_i)_{i=1}^k$ are independent, symmetric $\{-1,1\}$-valued random variables that are independent of $(X_i)_{i=1}^k$. If $T^\circ$ has a Rademacher type 2 constant $R_2(T^\circ)$, then
$$
\E \|k^{-1/2} \sum_{i=1}^k X_i\|_{T^\circ} =\E \|k^{-1/2} \sum_{i=1}^k \eps_i X_i\|_{T^\circ} \leq R_2(T^\circ) (\E \|X\|_{T^\circ}^2)^{1/2},
$$
and with probability at least $3/4$,
$$
{\rm diam}(T \cap {\rm ker}(\Gamma)) \leq \frac{c_3(\lambda)}{\sqrt{k}} \cdot \max\left\{ \E\|G\|_{T^\circ}, R_2(T^\circ) (\E \|X\|_{T^\circ}^2)^{1/2} \right\}.
$$
For example, if $T=B_1^n$ and $X \in \beta B_\infty^n$ almost surely, then $\|X\|_{T^\circ} = \|X\|_{\ell_\infty^n} \leq \beta$, $R_2(\ell_\infty^n) \leq \sqrt{\log n}$ and $\E\|G\|_{\ell_\infty^n} \lesssim \sqrt{\log n}$. Therefore,
$$
{\rm diam}(B_1^n \cap {\rm ker}(\Gamma)) \leq c_3(\lambda) \beta \sqrt{\frac{\log n}{k}}.
$$
\end{description}
\section{Proof of Theorem \ref{thm:main}} \label{sec:proof}
\begin{Lemma} \label{lemma:small-ball}
Let $\zeta$ be a random variable that satisfies
\begin{equation} \label{eq:single-small-ball}
Pr(|\zeta| \geq \lambda \|\zeta\|_{L_2}) \geq 1-\eps
\end{equation}
for constants $0<\eps<1/12$ and $\lambda>0$.

If $\zeta_1,...,\zeta_k$ are independent copies of $\zeta$, then with probability at least $1-2^{-6\eps k}$ there is a subset $J \subset \{1,...,k\}$ of cardinality at least $(1-6\eps)k$, and for every $j \in J$,
$$
|\zeta_j| \geq \lambda \|\zeta\|_{L_2}.
$$
\end{Lemma}

\proof It suffices to show that no more than $6\eps k$ of the $|\zeta_i|$'s are smaller than $\lambda \|\zeta\|_{L_2}$. By a binomial estimate, if $6\eps k \leq k/2$,
\begin{align*}
& Pr \left( \exists J \subset \{1,...,k\}, \ |J| = 6\eps k, \ |\zeta_j| \leq \lambda \|\zeta\|_{L_2} \ {\rm if } \ j \in J \right)
\\
\leq &\binom{k}{6\eps k} Pr^{6 \eps k} \left(|\zeta| \leq \lambda \|\zeta\|_{L_2} \right) \leq \left(\frac{e}{6 \eps}\right)^{6 \eps k} \cdot \eps^{6 \eps k} \leq 2^{-6 \eps k}.
\end{align*}
\endproof
Let $\{\zeta^i : 1 \leq i \leq N\}$ be a collection of random variables, and for every $i$ let $Z_i \in \R^k$ be a random vector with independent coordinates, distributed according to the random variable $\zeta^i$. Denote by $Z_i(j)$ the $j$-th coordinate of $Z_i$.
\begin{Corollary} \label{cor:large-set}
If each $\zeta^i$ satisfies the small-ball condition \eqref{eq:single-small-ball} and $N \leq 2^{3\eps k}$, then with probability at least $1-2^{-3\eps k}$, for every $1 \leq i \leq N$ there is a subset $J_i \subset \{1,...,k\}$, of cardinality at least $(1-6\eps)k$, and
$$
|Z_i(j)| \geq \lambda \|\zeta^i\|_{L_2} \ \ {\rm for \ every} \ j \in J_i.
$$
\end{Corollary}

\noindent{\bf Proof of Theorem \ref{thm:main}.} Let $\eps=1/600$ and observe that by the small ball assumption and since $X$ is isotropic,
$$
Pr( |\inr{X,x}| \geq \lambda\|x\|_{\ell_2^n} ) \geq 99/100 \geq (1-\eps)
$$
for every $x \in \R^n$.

Fix $r>0$ to be named later and set $T_r = T \cap r S^{n-1}$. Let
$$
\rho = c\frac{\E\|G\|_{T_r^\circ}}{\sqrt{\eps k}}
$$
for a suitable absolute constant $c$ and set $V_r \subset T_r$ to be a maximal $\rho$-separated subset of $T_r$ with respect to the $\ell_2^n$ norm.
Sudakov's inequality (see, e.g. \cite{Pis,LT}) shows that for the right choice of $c$, $|V_r| \leq 2^{3\eps k}$.

Let
$$
\zeta^i = \inr{X,v_i}, \ \ v_i \in V_r, \ \ 1 \leq i \leq 2^{3 \eps k},
$$
and set $Z_i=(\zeta_j^i)_{j=1}^k$, a vector whose coordinates are independent copies of $\zeta^i$.

Applying Corollary \ref{cor:large-set} to the set $\{Z_i: 1 \leq i \leq 2^{3\eps k}\}$, it follows that with probability at least $1-2^{-3\eps k}$, for every $v \in V_r$ there is a subset $J_v \subset \{1,...,k\}$, $|J_v| \geq (1-6\eps)k = 99k/100$, and for every $j \in J_v$,
\begin{equation} \label{eq:event-1}
|\inr{X_j,v}| \geq \lambda \|v\|_{\ell_2^n} = \lambda r,
\end{equation}
and the last equality holds because $V_r \subset r S^{n-1}$.

For every $x \in T_r$, let $\pi(x) $ be the nearest point to $x$ in $V_r$ with respect to the $\ell_2^n$ norm. Therefore,
$$
\E|\inr{X,x-\pi(x)}| \leq \|x-\pi(x)\|_{L_2(\mu)} = \|x-\pi(x)\|_{\ell_2^n} \leq \rho.
$$
By the Gin\'{e}-Zinn symmetrization inequality \cite{GZ84}, the contraction inequality for Bernoulli processes (see, e.g., \cite{LT}), and since $x-\pi(x) \in 2T \cap \rho B_2^n$ for every $x \in T_r$,
\begin{align*}
& \E \sup_{x \in T_r} \frac{1}{k} \sum_{i=1}^k |\inr{X_i,x-\pi(x)}|
\\
\leq & \rho + \E \sup_{x \in T_r} \left|\frac{1}{k} \sum_{i=1}^k |\inr{X_i,x-\pi(x)}|- \E|\inr{X_i,x-\pi(x)}| \right|
\\
\leq & \rho + \frac{2}{k} \E \sup_{x \in T_r} \left|\sum_{i=1}^k \eps_i |\inr{X_i,x-\pi(x)}|\right| \leq \rho + \frac{2}{k} \E \sup_{x \in T_r} \left|\sum_{i=1}^k \eps_i \inr{X_i,x-\pi(x)}\right|
\\
\leq & \rho + \frac{2}{k} \E \left\|\sum_{i=1}^k \eps_i X_i \right\|_{ (2T \cap \rho B_2^n)^\circ} \leq \rho + \frac{4}{k} \E \left\|\sum_{i=1}^k X_i \right\|_{ (T \cap \rho B_2^n)^\circ}.
\end{align*}
Hence, by the choice of $\rho$ and the trivial inclusion $T \cap \rho B_2^n \subset T$, \begin{equation} \label{eq:random-l-1}
\E \sup_{x \in T_r} \frac{1}{k} \sum_{i=1}^k |\inr{X_i,x-\pi(x)}| \leq c\frac{\E\|G\|_{T^\circ}}{\sqrt{\eps k}} + \frac{4}{k} \E \left\|\sum_{i=1}^k X_i \right\|_{T^\circ}.
\end{equation}
Set $A=\E \sup_{x \in T_r} k^{-1} \sum_{i=1}^k |\inr{X_i,x-\pi(x)}|$ and let
$$
W_r=\left\{ (\inr{X_i,x-\pi(x)})_{i=1}^k : x \in T_r\right\}.
$$
Note that for $0<\delta<1$, with probability at least $1-\delta$,
$$
W_r \subset k (A/\delta) B_1^k.
$$
On that event, if $(w_i)_{i=1}^k \in W_r$ and $(w_i^*)_{i=1}^k$ is a non-increasing rearrangement of $(|w_i|)_{i=1}^k$,
$$
w_{k/100}^* \leq \frac{\|w\|_{\ell_1^k}}{k/100}\leq \frac{100}{\delta}A.
$$
Thus, for every $x \in T_r$ there is a subset $J^\prime_x \subset \{1,...,k\}$ of cardinality at least $99k/100$, and for every $j \in J^\prime_x$,

\begin{equation} \label{eq:event-2}
|\inr{X_i,x-\pi(x)}| \leq \frac{100}{\delta} A.
\end{equation}

Fix $X_1,...,X_k$ in the intersection of the two events defined in \eqref{eq:event-1} and \eqref{eq:event-2}. For every $x \in T_r$ set $I_x=J^\prime_x \cap J_{\pi(x)}$. Observe that $|I_x| \geq 98k/100$ and that for every $i \in I_x$,
\begin{align*}
|\inr{X_i,x}| \geq & |\inr{X_i,\pi(x)}| - |\inr{X_i,x-\pi(x)}| \geq \lambda r - \frac{100}{\delta}A
\\
\geq & \lambda r - \frac{c_1}{\delta \sqrt{k}} \cdot \left(\E\|G\|_{T^\circ} +  \E \left\|\frac{1}{\sqrt{k}}\sum_{i=1}^k X_i \right\|_{T^\circ}\right).
\end{align*}
Therefore, if
$$
r \geq \frac{c_2(\lambda,\delta)}{\sqrt{k}} \cdot \max\left\{\E\|G\|_{T^\circ}, \E \left\|\frac{1}{\sqrt{k}}\sum_{i=1}^k X_i \right\|_{T^\circ}\right\},
$$
then with probability at least $1-\delta-2^{-3\eps k}= 1-\delta-2^{-k/200}$, for each $x \in T_r$, $|\inr{X_i,x}| \geq (\lambda/2) \|x\|_{\ell_2^n}$ on at least $98k/100$ coordinates; Thus,
\begin{equation} \label{eq:inproof-lower-bound}
\inf_{x \in T_r} \|\Gamma x\|_{\ell_2^k} \gtrsim \lambda \|x\|_{\ell_2^n}.
\end{equation}
Finally, using the convexity of $T$ and since the condition in \eqref{eq:inproof-lower-bound} is positive-homogeneous, \eqref{eq:inproof-lower-bound} holds for any $x \in T$ with $\|x\|_{\ell_2^n} \geq r$, as claimed.
\endproof
\section{concluding comments} \label{sec:comments}
The proof of Theorem \ref{thm:main} has two components. The first is based on a small-ball estimate for linear functionals and does not require additional information on their tails. Thus, this part holds even for heavy-tailed ensembles.

The more restrictive condition is on the random vector $k^{-1/2}\sum_{i=1}^k X_i$. Still, it is far easier to handle the norm $\|\sum_{i=1}^k X_i\|_{T^\circ}$ than the supremum of the quadratic empirical process indexed by $T$.

The estimate in Theorem \ref{thm:main} can be improved using what is, by now, a standard argument. First, observe that all the inequalities leading to \eqref{eq:random-l-1} hold in probability and not just in expectation (see, for example, \cite{VW,Dud-book}). Keeping the `localization' level $r$, one can define two fixed points:
\begin{equation*}
\rho_k(\delta,Q_1) = \inf \left\{ \rho:  Pr \left(\|k^{-1/2} \sum_{i=1}^k X_i \|_{ (T \cap \rho B_2^n)^\circ} \geq Q_1 \rho \sqrt{k} \right) \leq \delta\right\},
\end{equation*}
and
\begin{equation*}
r_k(Q_2) = \inf \{r : \E\|G\|_{(T \cap r S^{n-1})^\circ} \leq Q_2 r \sqrt{k}\}.
\end{equation*}
It is straightforward to verify that there are constants $Q_1$ and $Q_2$ that depend only on $\lambda$, for which, with probability at least $1-\delta-2^{-k/200}$, if
$$
\|x\|_{\ell_2^n} \gtrsim \max\{\rho_k(\delta,Q_1),r_k(Q_2)\},
$$
then
$$
\|\Gamma x\|_{\ell_2^k} \gtrsim \lambda \|x\|_{\ell_2^n}.
$$
Thus, on the same event,
$$
{\rm diam}(T \cap {\rm ker}(\Gamma)) \lesssim \max\{\rho_k(\delta,Q_1),r_k(Q_2)\}.
$$
Finally, it is possible to use a slightly more involved, empirical processes based method, that leads to an exponential probability estimate of $1-2\exp(-ck)$ in Theorem \ref{thm:main}. A result of a similar flavour, concerning the smallest singular value of a random matrix with iid rows may by found in \cite{KM}.

Since the goal in this note was to present the idea of using a simple small-ball argument, rather than pursuing an optimal result, we have opted to present this proof.

\end{document}